\documentclass[12pt]{amsart}
\usepackage{amsmath,amsthm,amsfonts,amssymb}
\begin{document} 
\newcommand{\B}{{\mathbb B}}
\newcommand{\C}{{\mathbb C}}
\newcommand{\N}{{\mathbb N}}
\newcommand{\Q}{{\mathbb Q}}
\newcommand{\Z}{{\mathbb Z}}
\renewcommand{\P}{{\mathbb P}}
\newcommand{\R}{{\mathbb R}}
\newcommand{\rc}{\subset}
\newcommand{\rank}{\mathop{rank}}
\newcommand{\trace}{\mathop{tr}}
\newcommand{\dimc}{\mathop{dim}_{\C}}
\newcommand{\Lie}{\mathop{Lie}}
\newcommand{\tensor}{\otimes}
\newcommand{\End}{\mathop{End}}
\newcommand{\Auto}{\mathop{{\rm Aut}_{\mathcal O}}}
\newcommand{\alg}[1]{{\mathbf #1}}
\newtheorem*{definition}{Definition}
\newtheorem*{claim}{Claim}
\newtheorem{corollary}{Corollary}
\newtheorem*{Conjecture}{Conjecture}
\newtheorem*{SpecAss}{Special Assumptions}
\newtheorem{example}{Example}
\newtheorem*{remark}{Remark}
\newtheorem*{observation}{Observation}
\newtheorem*{fact}{Fact}
\newtheorem*{remarks}{Remarks}
\newtheorem{lemma}{Lemma}
\newtheorem{proposition}{Proposition}
\newtheorem*{theorem}{Theorem}
\title{%
On a special class of complex tori
}
\author {J\"org Winkelmann}
\begin{abstract}
We investigate which complex tori admits complex Lie subgroups
whose closure is not complex.
\end{abstract}
%
\address{%
J\"org Winkelmann \\
 Institut Elie Cartan (Math\'ematiques)\\
 Universit\'e Henri Poincar\'e Nancy 1\\
 B.P. 239\\
 F-54506 Vand\oe uvre-les-Nancy Cedex\\
 France
}
\email{jwinkel@member.ams.org\newline\indent{\itshape Webpage: }%
http://www.math.unibas.ch/\~{ }winkel/
}
\maketitle
It is mentioned in an article of J.~Moser (\cite{M})
that for the torus
$T=(\C/\Z[i])^2$ every complex connected Lie subgroup of $T$ has
a closure which is a complex subtorus of $T$.
In general, i.e., if $T$ is an arbitrary compact complex
torus, the closure of a complex Lie subgroup in a compact
complex torus is a compact {\em real} subtorus which need not be
complex.
For instance, let
\[
T=(\C/\Z[i])\times (\C/\Z[\sqrt 2i])
\]
and let $H$ be the connected complex Lie subgroup which is the
image in $T$ of the diagonal line $\{(z,z):z\in\C\}$ in $\C^2$.
Then the preimage of $H$ in $\C^2$ can be described as
\[
\pi^{-1}(H)=
\{(z+n+mi,z+p+q\sqrt 2 i):z\in\C;n,m,p,q\in\Z\}
\]
whose closure is $\{(z,w)\in\C^2:\Re(z-w)=0\}$ and therefore of
real dimension three.
Thus $\bar H$ is a real subtorus of $T$ of real codimension one.

Our goal is to determine precisely the class of those compact complex
tori for which this phenomenon may occur.

\begin{theorem}
Let $T$ be a compact complex torus of dimension at least two.

Then the following two conditions are equivalent:

\begin{enumerate}
\item
For every connected complex Lie subgroup $H$ of $T$ the closure $\bar H$
in $T$ is a {\em complex} subtorus of $T$.
\item
There exists an elliptic curve $E$ with complex multiplication
 such that $T$ is isogenous to $E^n$ (with $n=\dim(T)$).
\end{enumerate}
\end{theorem}

We recall that an ``elliptic curve'' is a compact complex torus of
dimension one and that such an elliptic curve is said to have
``complex multiplication'' if $\End_\Q(E)$ is larger then $\Q$.

We recall that an endomorphism $f$ of an elliptic curve $E$ is a 
holomorphic Lie group homomorphism from $E$ to itself. Since $E$
is a commutative group, the set $\End(E)$ of all such endomorphisms
is a $\Z$-module in a natural way. Then $\End_\Q(E)$ is defined
as $\End(E)\tensor_\Z\Q$. Since every endomorphism $g$ of $E$ lifts
to an endomorphism of the universal covering $(\C,+)$, there is
a natural homomorphism from $\End_\Q(E)$ to $\C$.
Thus $\End_\Q(E)$ may be regarded as a $\Q$-sub algebra of $\C$.

For a given elliptic curve $E$ there are two possibilities:
Either $\End_\Q(E)=\Q$ or  $\End_\Q(E)$ is larger then $\Q$.
In the latter case $\End_\Q(E)$ 
(regarded as subalgebra of $\C$)
must contain a complex number $\lambda$ which
is not real. For this reason an elliptic curve is said to
have ``complex multiplication'' if $\End_\Q(E)\ne\Q$.

See e.g.~\cite{L} for more information about elliptic curves
with complex multiplication.

\begin{proof}
First let us assume property $(ii)$. Let $\lambda_0\in\End_\Q(E)$
with $\lambda_0\not\in\R$. 
Then $\lambda_0m\in\End(E)$ for some $m\in\N$.
Define $\lambda=\lambda_0m$. Let $E=\C/\Lambda$. Then we can realize
$T$ as a quotient $T=\C^n/\Gamma$ where $\Gamma$ is commensurable
with $\Lambda^n$. For every connected complex Lie subgroup $H\subset T$
we consider its preimage under the projection $\pi:\C^n\to T$.
Then $\pi^{-1}(H)=V+\Gamma$ for some complex vector subspace $V$
of $\C^n$. Now $\lambda\cdot(V+\Gamma)=
V+\lambda\cdot\Gamma$ and $V+\Gamma$ are commensurable.
Therefore the connected components of their respective closures in
$\C^n$ agree. It follows that the connected component $W$ of
the closure of $V+\Gamma$ in $\C^n$  is a real vector subspace
which is invariant under multiplication with $\lambda$.
Since $\lambda\in\C\setminus\R$, it follows that $W$ is a complex
vector subspace. Consequently $\Bar H=\pi(W)$ is a {\em complex}
Lie subgroup of $T$, i.e. a complex subtorus.

Now let us deal with the opposite direction.
As a preparation let us discuss real subtori of codimension one. 
If $S$ is a
real subtorus of $T$ of real codimension one, it corresponds
to a real hyperplane $H$ in $\C^n$. Then $H\cap iH$ is a complex
hyperplane in $\C^n$ which projects onto a connected complex Lie subgroup $A$
of $T$. By construction either this complex Lie subgroup $A$ is already closed
(i.e. a complex subtorus) or its closure equals $S$.

Let $U=\C^n$ and let $\P^*_\C(U)$ resp.~$\P^*_\R(U)$ denote the
spaces parametrizing the complex resp.\ real hyperplanes in $U$.
Then $H\mapsto H\cap iH$ defines a surjective continuous map
from $\P^*_\R(U)$ to $\P^*_\C(U)$. Let $\P^*_\Gamma(U)$ denote
the subset of those real hyperplanes which are generated by their
intersection with $\Gamma$. Observe that a $\R$-linear change of 
coordinates takes $\P^*_\Gamma(U)$ to $\P_{2n-1}(\Q)$ and
$\P_\R^*(U)$ to $\P_{2n-1}(\R)$. Therefore $\P^*_\Gamma(U)$
is dense in $\P_\R^*(U)$ and furthermore projects onto a dense
subset of $\P_\C^*(U)$.

Let us now assume condition $(i)$. Then for every real subtorus $S$ of
codimension one the connected complex Lie subgroup $A$ of codimension
one constructed above can not have $S$ as closure and therefore must
be complex compact subtorus. We thus obtain the following fact:

{\em Let $\P'$ denote the set of all complex hyperplanes in 
$\P_\C^*(U)$ which project onto compact complex subtori of $T$.
Then $\P'$ is dense in $\P_\C^*(U)$.}

As a consequence, there are compact complex subtori $(C_i)_{i=1..n}$
of codimension one such that the intersection $\cap_i C_i$
is discrete. It follows that there is a surjective homomorphism
of complex tori with finite kernel
\[
\psi: T \to \Pi_{i=1}^n \left(T/C_i\right).
\]
Thus $T$ is isogenous to a product of elliptic curves.

Let us now discuss two-dimensional quotient tori of $B$.
If $\tau:T\to B$ is a projection onto a two-dimensional torus
and $L\subset B$ is a one-dimensional complex Lie subgroup,
then $\overline{\tau^{-1}(L)}=\tau^{-1}(\bar L)$. Thus assuming
condition $(i)$ for $T$ implies the same condition for $B$.
If we now define $\P'_B$ as the subset of those complex lines
in $\C^2$ whose image in $B$ are complex subtori, then we obtain
that $\P'_B$ must be dense in $\P_1(\C)$.

Thus let us discuss $\P'_B$ for $B=E'\times E''$ where
$E'$ and $E''$ are elliptic curves.  
If $E'$ is not isogenous to $E''$, then $E'\times\{e\}$ and
$\{e\}\times E''$ are the only subtori of $B$ and $\P'_B$
can not be dense. Thus we may assume that $E'$ is isogenous to
$E''$.

If $E'$ does not have complex multiplication, then
$\P'_N\simeq\P_1(\Q)$ whose closure is $\P_1(\R)$ and which
therefore is not dense.

This leaves the case where $E'$ is isogenous to $E''$ and 
has complex multiplication.

If this is to hold for any two-dimensional quotient torus of $T$,
it requires that all the $E_i$ are isogenous to each other and
have complex multiplication.
\end{proof}


\begin{thebibliography}{Bla}

\bibitem{M}
Moser, J.:
On the persistence of pseudo-holomorphic curves on an 
almost complex torus (with an appendix by J\"urgen P\"oschel). 
\sl Invent. Math. \bf 119 \rm (1995), no. 3, 401--442.

\bibitem{L}
Lang, S.:
Complex multiplication.
GTM {\bf 255}. Springer. 1983.
\end{thebibliography}
\end{document}